# Some Methods to Solve the Diophantine Equation $ay^2 = b(x^k - c^k)$

XIAODONG ZHUANG[1], NIKOS E. MASTORAKIS[2]

[1] Qingdao University, Qingdao, CHINA
[2] Technical University of Sofia, Sofia, BULGARIA

*Abstract:* - The solution of the Diophantine equation $ay^2 = b(x^k - c^k)$ is studied in different cases of *k*. Methods are proposed for different cases of *k*, which are convenient for computer algorithm implementation. The original equation is reduced to lower-degree equarions, which can improve the efficiency in searching possible integer solutions. A sufficient condition of the unsolvability of $ay^2 = b(x^{2n} - c^{2n})$ is presented.

## 1 Introduction

The Hilbert's tenth problem asked whether there is a general algorithm to determine the solvability of any Diophantine equation. A negative answer has been proved [1][2]. Due to the importance of Diophantine equation in theory and application [3][4][5][6][7], it is still meaningful to study the solution of specific types Diophantine equations, and it is also possible to find effective methods to solve specific type of them. Compared to the equations with concrete coefficients, the work on those with undetermined coefficients is more general. In this paper, the type of equation $ay^2 = b(x^k - c^k)$ is studied (where *a*, *b* and *k* are positive integers.)

Baker's method (and later improvements) may estimate the bound of possible integer solutions for Diophantine equations, while the range estimated is often very large [8][9]. To verify each integer in the range by enumeration, the assistance of computer programming is necessary, but there is still huge amount of computation. In this paper, the equation $ay^2 = b(x^k - c^k)$ is studied by factorization of the polynomial on the right side, and effective methods are presented to improve the efficiency in finding possible solution. Moreover, a sufficient condition of the unsolvability of $ay^2 = b(x^{2n} - c^{2n})$ is presented.

## 2 The Solution of $y^2 = x^{2n} - c^{2n}$

In this section, $y^2 = x^{2n} - c^{2n}$ is solved in finite computation steps (without the pre-estimation of the range of integer solution). Factorize the right side of this equation:

$$y^2 = x^{2n} - c^{2n} = (x^n + c^n)(x^n - c^n) \quad (1)$$

Suppose $d = \gcd(x^n + c^n, x^n - c^n)$. We have: $d \mid [(x^n + c^n) - (x^n - c^n)]$ and $d \mid (2c^n)$. Therefore, all the possible values of *d* belongs to the factor set of $2c^n$, which can be enumerated finitely. Transform (1) to seperate the coprime parts of the right side:

$$\begin{aligned} y^2 &= (x^n + c^n)(x^n - c^n) \\ &= d^2 \cdot \frac{x^n+c^n}{d} \cdot \frac{x^n-c^n}{d} \end{aligned} \quad (2)$$

It follows that

$$\left(\frac{y}{d}\right)^2 = \frac{x^n+c^n}{d} \cdot \frac{x^n-c^n}{d} \quad (3)$$

where the left side is a square number, and the right side has two coprime items. Therefore, both items of the right side are square umbers respectively:

$$\begin{cases} \frac{x^n+c^n}{d} = u^2 \\ \frac{x^n-c^n}{d} = v^2 \end{cases} \quad u, v \in Z \quad (4)$$

and

$$\begin{cases} x^n + c^n = du^2 \\ x^n - c^n = dv^2 \end{cases} \quad (5)$$

Eliminate *x* in (5), we have $du^2 - dv^2 = 2c^n$ and

$$u^2 - v^2 = \frac{2c^n}{d} \quad (6)$$

which can be rewritten as:

$$(u + v)(u - v) = \frac{2c^n}{d} \quad (7)$$

The right side of (7) can be factorizd as the multiplication of two facors: $\frac{2c^n}{d} = g \cdot h$, and the left side is also a multiplication of two items. Therefore, the solution can be searched by enumeration of all possible combinations of factor correspondence between the left and right side:

$$\begin{cases} u + v = g \\ u - v = h \end{cases} \quad (8)$$

It follows that $u = \frac{g+h}{2}$, $v = \frac{g-h}{2}$. If $\frac{g+h}{2}$ and $\frac{g-h}{2}$ are both integers, check whether $\sqrt[n]{d \cdot u^2 - c^n}$ is an integer. If so, one solution of the original equation is found:

$$x = \pm\sqrt[n]{d \cdot u^2 - c^n}, y = \pm\sqrt{x^{2n} - c^{2n}} \quad (9)$$

Therefore, by enumerating $d$ (i.e. the factor of $2c^n$) and $(g,h)$ (i.e. $\frac{2c^n}{d} = g \cdot h$), $y^2 = x^{2n} - c^{2n}$ can be solved, conveniently by computer programming.

## 3 The Solution of $ay^2 = b(x^{2n} - c^{2n})$

The method in the previous section can be extended to the equation with arbitrary integer coefficients:

$$ay^2 = b(x^{2n} - c^{2n}) \quad (10)$$

Make the left side a square:

$$(ay)^2 = ab(x^n + c^n)(x^n - c^n) \quad (11)$$

and the coprime parts of the right side is seperated as following.

Suppose the maximum square factor of $ab$ is $r^2$. We have $ab = r^2 t$, where $t$ is the square-free factor which equals the multiplication of distinct prime numbers.

Suppose $d = \gcd(x^n + c^n, x^n - c^n)$, we have $d|2c^n$. Therefore, $d$ belongs to the factor set of $2c^n$, which can be finitely enumerated. Then we have:

$$\begin{aligned}(ay)^2 &= ab(x^n + c^n)(x^n - c^n) \\ &= (r^2 t) \cdot d^2 \cdot \frac{x^n + c^n}{d} \cdot \frac{x^n - c^n}{d} \\ &= (rd)^2 \cdot t \cdot \frac{x^n + c^n}{d} \cdot \frac{x^n - c^n}{d}\end{aligned} \quad (12)$$

and

$$\left(\frac{ay}{rd}\right)^2 = t \cdot \left(\frac{x^n + c^n}{d}\right) \cdot \left(\frac{x^n - c^n}{d}\right) \quad (13)$$

The left side of (13) is a square number. On the right side, the two items in brackets are coprime, and $t$ is square-free factor. Therefore, we have

$$\left(\frac{x^n + c^n}{d}\right) \cdot \left(\frac{x^n - c^n}{d}\right) = u^2 v^2 t \quad (14)$$

Otherwise, the right side of (13) can not equal a square number. Again because the above two items in brackets are coprime, they must have the following form:

$$\begin{cases}\frac{x^n + c^n}{d} = gu^2 \\ \frac{x^n - c^n}{d} = hv^2\end{cases} \quad (15)$$

where $t=gh$. Because $t$ is square-free factor, $\gcd(g,h)=1$. And all possible $(g,h)$ can be enumerated finitely. Then we have

$$\begin{cases}x^n + c^n = dgu^2 \\ x^n - c^n = dhv^2\end{cases} \quad (16)$$

Eliminate $x$ in (16), we have $dgu^2 - dhv^2 = 2c^n$. Because $d|2c^n$, we have $gu^2 - hv^2 = \frac{2c^n}{d}$. Make the first item on the left side a square: $(gu)^2 - ghv^2 = \frac{2gc^n}{d}$, and substitue $t=gh$ into it:

$$(gu)^2 - tv^2 = \frac{2gc^n}{d} \quad (17)$$

Let $w=gu$, (17) is rewritten as:

$$w^2 - tv^2 = \frac{2gc^n}{d} \quad (18)$$

which is a generalized Pell equation (with the restriction $g|w$). To solve (18), the methods already proposed for generalized Pell equations is ready to apply here [10][11][12][13][14]. If (18) has no integer solution, the original equation also has no integer solution. If (18) has interger solution, it has infinite solutions, based on which the solution of the original equation can be searched. Here the possible range of solution to (18) can be estimated (by Baker methods and the later improvements) [8][9], and in this range check the following to determine the solution of (10):

(1) whether $g|w$ holds true;
(2) whether $u = \frac{w}{g}$ and $v$ are coprime;
(3) whether $\sqrt[n]{dgu^2 - c^n}$ is an integer.

If all the above three are satisfied, a solution of (10) is found:

$$x = \pm\sqrt[n]{dgu^2 - c^n},\ y = \pm\sqrt{\frac{b}{a}(x^{2n} - c^{2n})} \quad (19)$$

Therefore, by enumerating $d$ (i.e. the factor of $2c^n$) and $(g,h)$ (i.e. $t=gh$ and $t$ is the square-free factor of $ab$), the solution of (10) can be searched by computer programming. Moreover, based on the above analysis, a sufficient condition of the unsolvability of (10) is presented as the following:

*For $ay^2 = b(x^{2n} - c^{2n})$, let $t$ be the square-free factor of ab. If for any t's factor g, and any ($2c^n$)'s factor d, the equation $w^2 - tv^2 = \frac{g \cdot 2c^n}{d}$ has no integer solution, then $ay^2 = b(x^{2n} - c^{2n})$ has no integer solution either.*

Because when $d$ traverses all the factors of $2c^n$, $\frac{2c^n}{d}$ also traverses all the factors of $2c^n$, the sufficient condition of (10)'s unsolvability can be rewritten as:

*For $ay^2 = b(x^{2n} - c^{2n})$, let $t$ be the square-free factor of ab. If for any t's factor g, and any ($2c^n$)'s factor e, the equation $w^2 - tv^2 = ge$ has no integer solution, then $ay^2 = b(x^{2n} - c^{2n})$ has no integer solution either.*

## 4 The Solution of $ay^2 = b(x^{2n+1} - c^{2n+1})$

Facorize $(x^{2n+1} - c^{2n+1})$ on the right side, and make the left side a square number, we have

$$\begin{aligned}(ay)^2 &= ab(x - c)(x^{2n} + cx^{2n-1} + \cdots + c^{2n}) \\ &= ab(x - c) \cdot \sum_{i=0}^{2n} c^i x^{2n-i}\end{aligned} \quad (20)$$

Let $v=x-c$, (20) is rewritten as:

$$(ay)^2 = abv \cdot \sum_{i=0}^{2n} c^i (v+c)^{2n-i} \qquad (21)$$

Then expand the summation item on the right side:

$$\sum_{i=0}^{2n} c^i (v+c)^{2n-i}$$
$$= \sum_{i=0}^{2n} c^i \left(\sum_{j=0}^{2n-i} \boldsymbol{C}_{2n-i}^{j} \cdot v^j \cdot c^{2n-j}\right)$$
$$= \sum_{m=0}^{2n} \left(\sum_{s=m}^{2n} \boldsymbol{C}_s^m \cdot c^{2n-m}\right) v^m$$

Note that the constant item in the above summation is $(2n+1)c^{2n}$. Let $q_m = \sum_{s=m}^{2n} \boldsymbol{C}_s^m \cdot c^{2n-m}$, (21) is rewritten as:

$$(ay)^2 = ab \cdot v \cdot \sum_{m=0}^{2n} q_m v^m \qquad (22)$$

Suppose the maximum square factor of $ab$ is $r^2$. Then $ab = r^2 t$, where $t$ is the square-free factor. Suppose $d = \gcd(v, \sum_{m=0}^{2n} q_m v^m)$, then it follows that $d \mid ((2n+1)c^{2n})$. Note that $(2n+1)c^{2n}$ is the constant item in $\sum_{m=0}^{2n} q_m v^m$. Then all the possible values of $d$ belongs to the factor set of $(2n+1)c^{2n}$. Therefore, (22) is rewritten as:

$$(ay)^2 = (r^2 t) \cdot d^2 \cdot \frac{v}{d} \cdot \left(\frac{1}{d}\sum_{m=0}^{2n} q_m v^m\right)$$

or

$$\left(\frac{ay}{rd}\right)^2 = t \cdot \left(\frac{v}{d}\right) \cdot \left(\frac{1}{d}\sum_{m=0}^{2n} q_m v^m\right) \qquad (23)$$

where the left side is a square number, and the two items in brackets on the right side are coprime. And $t$ is a square-free factor, so the following must holds true:

$$\left(\frac{v}{d}\right) \cdot \left(\frac{1}{d}\sum_{m=0}^{2n} q_m v^m\right) = e^2 f^2 t \qquad (24)$$

Otherwise the right side of (23) can not be a square number. Because $t$ can be factorized as the multiplication of two coprime factors $t=gh$ (all the possible factorization can be enumerated), the two items in brackets in (24) must be:

$$\begin{cases} \dfrac{v}{d} = ge^2 \\ \dfrac{1}{d}\displaystyle\sum_{m=0}^{2n} q_m v^m = hf^2 \end{cases} \qquad e, f \in Z$$

and we have

$$\begin{cases} v = dge^2 \\ \sum_{m=0}^{2n} q_m v^m = dhf^2 \end{cases} \qquad (25)$$

Based on (25), a method to search the solutions is presented. First, estimate the range of solution by Baker's method (and the later improvements) [8][9], which also determins the range of $e$. Then traverse this range of $e$, calculate $v=dge^2$, substitute $v$ into the second equation of (25), and check whether the second equation holds (i.e whether $\frac{\sum_{m=0}^{2n} q_m v^m}{dh}$ is a square number). If so, we get one solution of $ay^2 = b(x^{2n+1} - c^{2n+1})$:

$$x = v + c,\ y = \pm\sqrt{\frac{b}{a}(x^{2n+1} - c^{2n+1})}$$

Therefore, by enumerating $d$ (i.e. the factor of $(2n+1)c^{2n}$) and $(g,h)$ (i.e. $t=gh$), the solution of $ay^2 = b(x^{2n+1} - c^{2n+1})$ can be searched by computer programming.

# 5 Conclusion

In this paper, factorization and GCD (greatest common factor) enumeration for integer polynomials is utilized to solve the Diophantine equation $ay^2 = b(x^k - c^k)$. The original equation is converted to generalized Pell equation or lower degree equation group, which can be solved by existing methods in a more easy way. Effective algorithms can be implemented by computer programming based on the methods prposed, by which the integer solutions of the original equation can be searched more efficiently. At the same time, a sufficient consition of the unsolvability of $ay^2 = b(x^k - c^k)$ is prposed.

The method proposed can be easily extended to the equation $ay^{2m} = b(x^k - c^k)$. Future study may improve the sufficient conditions, and introduce detailed congruence analysis into the methods to further improve the efficiency.